\documentclass[letterpaper, 10 pt, conference]{ieeeconf}

\IEEEoverridecommandlockouts                              
\usepackage{tikz}
\usepackage{subcaption}
\usepackage{booktabs}
\usepackage{cite}

\usepackage{algorithm}
\usepackage{makecell}
\usepackage{algpseudocode}
\usepackage{amsmath,amssymb,amsfonts}
\usepackage{graphicx}
\usepackage{textcomp}
\def\BibTeX{{\rm B\kern-.05em{\sc i\kern-.025em b}\kern-.08em
    T\kern-.1667em\lower.7ex\hbox{E}\kern-.125emX}}
\usepackage{mathrsfs}
\usepackage{dsfont}
\newtheorem{thm}{\bf Theorem}[section]
\newtheorem{cor}[thm]{\bf Corollary}

\newtheorem{rem}[thm]{\bf Remark}

\newtheorem{prop}[thm]{\bf Proposition}

\newtheorem{deff}[thm]{\bf Definition}
\newtheorem{eg}[thm]{\bf Example}

\renewcommand{\ll}{\mathcal{L}}
\renewcommand{\lll}{\mathcal{L}_{\lambda}}
\newcommand{\X}{\mathcal{X}}

\renewcommand{\S}{\mathcal{S}}
\newcommand{\I}{\mathcal{I}}

\newcommand{\B}{\mathcal{B}}
\newcommand{\rr}{\mathcal{R}}
\renewcommand{\S}{\mathcal{S}}

\newcommand{\K}{\mathcal{K}}
\newcommand{\T}{\mathcal{T}}

\newcommand{\F}{\mathcal{F}}
\newcommand{\R}{\mathbb{R}}
\newcommand{\N}{\mathbb{N}}
\newcommand{\Z}{\mathbb{Z}}
\newcommand{\C}{\mathbb{C}}

\newcommand{\Tb}{\mathbf{T}}
\newcommand{\Lb}{\mathbf{L}}

\newcommand{\tomg}{\tilde{\omega}}
\newcommand{\xb}{\mathbf{x}}
\newcommand{\vb}{\mathbf{\xi}}

\newcommand{\Zk}{\mathcal{Z}}
\newcommand{\zk}{\mathfrak{z}}

\newcommand{\eps}{\varepsilon}
\newcommand{\tta}{\vartheta}

\newcommand{\Nb}{\mathbf{N}}

\newcommand{\slim}{\mathop{\mathrm{s\text{-}\lim}}}
\newcommand{\ra}{\rightarrow}
\newcommand{\sra}{\rightharpoonup}

\newcommand{\id}{\operatorname{I}}
\newcommand{\dom}{\operatorname{dom}}
\newcommand{\inte}{\operatorname{int}}

\renewcommand{\Re}{\operatorname{Re}}
\newcommand{\set}[1]{\left\{#1\right\}}

\newcommand{\norm}[1]{\left\Vert #1 \right\Vert}
\newcommand{\abs}[1]{\left\vert #1 \right\vert}

\newcommand{\BlackBox}{\rule{1.5ex}{1.5ex}}    
\newcommand{\pfbox}{\hfill\BlackBox}    
\newcommand{\Qed}{\hfill$\diamond$}

\newenvironment{pf}{\par\noindent{\bf Proof:\ }}{\hfill\BlackBox\\[2mm]}

\usepackage{hyperref}
 \hypersetup{
    colorlinks=true,
    linkcolor=blue,
    filecolor=magenta,      
    urlcolor=cyan,
    pdftitle={Overleaf Example},
    pdfpagemode=FullScreen,
    }
    \usepackage[misc]{ifsym}

\title{\bfseries Koopman-Based Learning of Infinitesimal Generators \\without Operator Logarithm}

\author{Yiming Meng, Ruikun Zhou, Melkior Ornik, and Jun Liu 
\thanks{
This research was supported in part by NASA under grant numbers 80NSSC21K1030 and 80NSSC22M0070, by the Air Force Office of Scientific Research under grant number FA9550-23-1-0131, and in part by the Natural Sciences and Engineering Research Council of Canada and the Canada Research Chairs program.}
\thanks{
Yiming Meng is with  the
Coordinated Science Laboratory, University of Illinois Urbana-Champaign,
Urbana, IL 61801, USA.
        {\tt\small ymmeng@illinois.edu}.
}
\thanks{
Ruikun Zhou is with the Department of Applied Mathematics,
University of Waterloo, Waterloo ON N2L 3G1, Canada {\tt\small ruikun.zhou@uwaterloo.ca}.}
\thanks{Melkior Ornik is with the Department of Aerospace Engineering and the
Coordinated Science Laboratory, University of Illinois Urbana-Champaign,
Urbana, IL 61801, USA.
        {\tt\small mornik@illinois.edu}.}
\thanks{
Jun Liu is with the Department of Applied Mathematics,
University of Waterloo, Waterloo ON N2L 3G1, Canada {\tt\small j.liu@uwaterloo.ca}.}
}

\begin{document}

\maketitle

\begin{abstract}
To retrieve transient transition information of unknown systems from discrete-time observations, the Koopman operator structure has gained significant attention in recent years, particularly for its ability to avoid time derivatives through  the Koopman operator logarithm. However, the effectiveness of these logarithm-based methods has only been demonstrated within a restrictive function space. In this paper, we propose a logarithm-free technique for learning the infinitesimal generator without disrupting the Koopman operator learning framework.

\end{abstract}

\begin{keywords}
Unknown nonlinear systems, Koopman operators, infinitesimal generator, 
system identification, verification.
\end{keywords}

\section{Introduction}
\label{sec:introduction}

Verification of dynamical system properties and achieving autonomy are two important directions for the future of industrial intelligence, with applications in numerous fields, including mathematical finance, automated vehicles, power systems, and other physical sciences. 

Witnessing the success of problem-solving within the data paradigm, there has been a surge of interest in revealing the governing equations of continuous-time dynamical systems from time-series data to better understand the underlying physical laws \cite{sjoberg1995nonlinear, haber1990structure}.  Additional interests in safety-critical industries include data-driven stability and safety analysis, prediction, and control. Techniques such as Lyapunov and barrier certificates have proven effective \cite{mauroy2013spectral, mauroy2016global, deka2022koopman, proctor2016dynamic, meng2023learning}. It is worth noting, however, that these practical concerns require information on the vector fields, the value functions that abstract system performance, and the corresponding Lie derivatives, all underpinned by an understanding of the infinitesimal generator \cite{mitchell2005time, lin1996smooth, liu2023physics, meng2022smooth, meng2023lyapunov}. Considering nonlinear effects, challenges therefore arise in the converse identification of infinitesimal system transitions based on discrete-time observations that represent cumulative trajectory behaviors.

Direct methods, such as the sparse identification of nonlinear dynamics (SINDy) algorithm \cite{brunton2016discovering}, have been developed to identify state dynamics by relying on nonlinear parameter estimation \cite{nash1987nonlinear, varah1982spline} and static linear regression techniques. However, the accurate approximation of time derivatives of the state may not be achieved due to challenges such as low sampling rates, noisy measurements, and short observation periods. 
Furthermore, the data cannot be reused in the proposed structure for constructing other value functions (e.g., Lyapunov or barrier functions) for stability, reachability, and safety analysis. This limitation extends to verifying their Lie derivatives along the trajectories, which is crucial for demonstrating the evolving trends of the phase portraits.

Comparatively, the operator logarithm-based Koopman generator learning structure \cite{mauroy2019koopman, klus2020data,drmavc2021identification, black2023safe} does not require the estimation of time derivatives, enabling a data-driven estimation of Lie derivatives. This approach can potentially circumvent the need for high sampling rates and extended observation periods. Heuristically, researchers tend to represent the Koopman operator $\K_t$ by an exponential form of its infinitesimal generator $\ll$ as $\K_t=e^{t\ll}$, leading to the converse representation $\ll = \frac{1}{t}\log(\K_t)$ for any $t>0$. However, representing Koopman operators in exponential form requires the boundedness of the generator. Additionally, the operator logarithm is a single-valued mapping only within a specific sector of the spectrum.
Recent studies \cite{zeng2022sampling, zeng2023generalized} have investigated the sufficient and necessary conditions under which the Koopman-logarithm-based generator learning method can be uniquely identified. However, these conditions are less likely to be verifiable for unknown systems.

In this paper, we introduce a logarithm-free generator learning scheme that does not require knowledge of the spectrum properties and present preliminary results.   
This scheme will be compatible with the current advances in \cite{meng2023learning} for Koopman-based construction of maximal Lyapunov functions. It is important to note that the method in \cite{meng2023learning} assumes full knowledge of the equilibrium point and acknowledges that verification of the constructed Lyapunov function depends on information about the actual system transitions, which may diminish its predictive value in stability analysis. 

The rest of the paper is organized as follows. Preliminaries are covered in Section \ref{sec: pre}. The formal finite-horizon approximation of the generator using Koopman operators is presented in Section \ref{sec: characterize}. The data-driven algorithm is discussed in Section \ref{sec: alg}. Finally, case studies are presented in Section \ref{sec: num}. 

Due to page limitations and to enhance readability, this paper presents only the primary results of the proposed method. Detailed proofs are omitted and will be published elsewhere.

\textbf{Notation:} 
We denote by $\R^n$ the Euclidean space of dimension $n>1$, and by $\R$ the set of real numbers. 
For $x\in\R^n$ and $r\ge 0$, we denote the ball of radius $r$ centered at $x$ by $\B(x, r)=\set{y\in\R^n:\,\abs{y-x}\le r}$, where $\abs{\;\cdot\;}$ is the Euclidean norm. For a closed set $A\subset\R^n$ and $x\in\R^n$, we denote the distance from $x$ to $A$ by $\abs{x}_{A}=\inf_{y\in A}\abs{x-y}$ and $r$-neighborhood of $A$ by $\B(A, r)=\cup_{x\in A}\B(x, r)=\set{x\in\R^n:\,\abs{x}_A\le r}$. For a set $A\subseteq\R^n$, $\overline{A}$ denotes its closure, 
$\inte(A)$ denotes its interior  and $\partial A$ denotes its boundary.  
For finite-dimensional matrices, we use the Frobenius norm  $\|\cdot\|_F$ as the metric.  
Let $C(\Omega)$ be the set of continuous functions with domain $\Omega$. We denote the set of $i$ times continuously
differentiable functions by $C^i(\Omega)$. 



\section{Preliminaries}\label{sec: pre}
\subsection{Dynamical Systems}
Given a pre-compact state space $\X\subseteq\R^n$, 
we consider a continuous-time nonlinear dynamical system of the form 

\begin{equation}\label{E: sys}
    \dot{\xb}(t) = f(\xb(t)),\;\;\xb(0)=x\in\X,\;t\in[0,\infty), 
\end{equation}
where $x$ denotes the initial condition, and the vector field $f:\X\ra\X$ is assumed to be locally Lipschitz continuous. We denote by $\phi: \I\times \X\rightarrow \X$ the forward flow, also known as the solution map.

The evolution of observable functions of system \eqref{E: sys} restricted to $C(\X)$ is governed by the family of Koopman operators, defined by \begin{align}
\mathcal{K}_t h = h\circ \phi(t, \cdot), \quad h \in C(\X)
\end{align}
for each $t\geq 0$, where $\circ$ is the composition operator. By the properties of the flow map, it can be easily verified that $\set{\K_t}_{t\geq 0}$ forms a semigroup, with its (infinitesimal) generator defined by  
\begin{equation}\label{E: generator}
    \ll h(x):= \lim_{t\ra 0}\frac{\K_t h(x)-h(x)}{t}. 
\end{equation}
Within the space of continuous observable functions, the limit in Eq. \eqref{E: generator} exists on the subspace of all continuously differentiable functions, indicating that $\dom(\ll)=C^1(\X)$. In this case, we can verify that the generator of Koopman operators is such that $$\ll h(x) = \nabla h(x) \cdot f(x). $$ 

Note that the finite-difference method \cite{bramburger2024auxiliary, nejati2021data, wang2022data} indirectly approximates the generator by setting a small terminal   time without taking the limit on the r.h.s. of \eqref{E: generator}. Through this approximation scheme of the time derivative, it can be anticipated that the precision heavily depends on the size of that terminal time.


\subsection{Representation of Koopman Operators} 
If $\ll$ is a
bounded linear operator that generates $\{\K_t\}$, then $\K_t = e^{t\ll}$ for each $t$ in the uniform (operator) topology.  However, many of the differential operators that occur in mathematical physics are unbounded. We revisit some concepts to show how to interpret $\K_t$ in terms of the possibly unbounded generator on an exponential scale. 

\begin{deff}[Resolvents]
The resolvent set of $\ll$ is defined as 
$     \rho(\ll):=\set{\lambda\in\C: \lambda\id-\ll\;\text{is invertible}}$. 
For any $\lambda\in\rho(\ll)$, the resolvent  operator is defined as
\begin{equation}
    R(\lambda; \ll):=(\lambda \id-\ll)^{-1}, 
\end{equation}
which is a bounded linear operator \cite[Chap. I,  Theorem 4.3]{pazy2012semigroups}.
\Qed
\end{deff}

We further define the \textit{Yosida approximation} of  $\ll$  as
\begin{equation}
    \ll_\lambda:=\lambda \ll R(\lambda;\ll) = \lambda^2R(\lambda;\ll)-\lambda \id. 
\end{equation}
Note that $\set{\ll_\lambda}_{\lambda\in\rho(\ll)}$ is a family of bounded linear operators, and $e^{t\ll_\lambda}$ is well-defined for each $\lambda\in\rho(\ll)$. 

\begin{deff}[Strong convergence]Let $(\F, \|\cdot\|_\F)$ be a Banach space. 
    Let $B: \F\ra\F$ and $B_n:\F\ra\F$, for each $n\in\N$, 
be   linear operators. Then, the $\{B_n\}_{n\in\N}$ is said to converge to $B$ \emph{strongly}, denoted by $\B_n\sra B$, if $\lim_{n\ra\infty}\|B_nh-Bh\|_\F=0$ for each $h\in\F$. We also write $B=\slim_{n\ra\infty}B_n$. \Qed
    
\end{deff}

\begin{thm}\label{thm: H-Y} \cite[Chap. I, Theorem 5.5]{pazy2012semigroups}
$\K_t  = \slim_{\lambda\ra\infty} e^{t\ll_\lambda}$ for all $t\geq 0$ on the uniform topology. 
\end{thm}


\section{The Formal Approximation of the Infinitesimal
Generators}\label{sec: characterize}
Section \ref{sec: pre} presents preliminary results on representing Koopman operators using the asymptotic approximation of $e^{t\ll_\lambda}$ as $\lambda$ approaches infinity. In this section, we aim to find the correct converse representation of $\ll$ based on $\{\K_t\}_{t\geq 0}$.  The ultimate goal is to leverage the Koopman learning structure to learn the generator. 

Note that it is intuitive to take the operator logarithm such that $\ll = \frac{1}{t}\log \K_t$ when $\ll$ is bounded.  The spectrum's sector should be confined to make the logarithm a single-valued mapping \cite{zeng2022sampling, zeng2023generalized}. However, for an unbounded $\ll$, there is no direct connection.  In this section, we establish the connection for how the unbounded $\ll$ can be properly approximated based on $\{\K_t\}_{t\geq 0}$.

\subsection{Asymptotic Approximations of Generators}\label{sec: app_generator}
Let $C(\X)$ be endowed with the uniform norm $|\cdot|_\infty:=\sup_{x\in\X}|\cdot|$. 
Suppose that $\{\K_t\}_{t\geq 0}$ is a contraction semigroup where $$ \|\K_t\|:=\sup_{\|h\|_{\F}=1}|\K_th|_\infty\leq 1$$  for all $t\geq 0$. Then, $\rho(\ll)\supseteq (0, \infty)$, and by \cite[Lemma 1.3.3]{pazy2012semigroups}, we have that $\set{\lll}_{\lambda>0}$ converges strongly to $\ll$. Specifically, when working on the $C^2(\X)$ subspace, the convergence rate is given by $\mathcal{O}(\lambda^{-1})$. 

To address more general cases, we first present the following facts. 

\begin{prop}\label{prop: fact1}
   For system \eqref{E: sys}, there exist  constants $\omega\geq 0$ and $M\geq 1$ such that $\norm{\K_t}\leq Me^{\omega t}$ for all $t\geq 0$. In addition, for any $\lambda\in\C$, the family $\set{e^{\lambda t}\K_t}_{t\geq 0}$ 
   is a $C_0$ semigroup with generator $\ll+\lambda \id: \dom(\ll)\ra C(\X)$. 
\end{prop}

Intuitively, $M$ represents the uniform scaling of the magnitude of the Koopman operator, while $\omega$  indicates the dominant exponential growth or decay rate of the flow on $C(\X)$.  The above properties provide us with a tool to convert the dominant exponential growth or decay rate of the flow on $C(\X)$ by introducing an extra term $\lambda \id$ to the original generator. Specifically, the transformation $e^{-\omega t}\K_t$ allows us to work in a topology (by shifting the direction of flows and uniformly compressing by $M$) that is equivalent to the uniform topology of continuous functions, where the semigroup $\{e^{-\omega t}\K_t\}$ is a contraction. The convergence of $\slim_{\lambda\ra\infty}\lll = \ll$ for $\set{\lll}_{\lambda>\omega}$ and its reciprocal convergence rate can similarly be demonstrated in this new topology.

\subsection{Representation of  Resolvent Operators}

Motivated by representing $\ll$ by $\set{\K_t}_{t\geq 0}$ and  the Yosida approximation for $\ll$ on $\{\lambda>\omega\}$, we  establish a connection between $R(\lambda;\ll)$ and $\set{\K_t}_{t\geq 0}$. 

\begin{prop}\label{prop: R_form}
    Let $R(\lambda)$ on $C(\X)$ be defined by
    \begin{equation}
        R(\lambda)h:=\int_0^\infty e^{-\lambda t}(\K_t h) dt. 
    \end{equation}
    Then, for all $\lambda>\omega$,  
    \begin{enumerate}
        \item $R(\lambda)(\lambda \id-\ll)h = h$ for all $h\in C^1(\X)$;
        \item $(\lambda \id-\ll)R(\lambda)h = h$ for all $h\in C(\X)$.
    \end{enumerate}
\end{prop}

Restricted to the $C^1(\X)$ subspace, to use the approximation in Section \ref{sec: app_generator}, we can replace  $R(\lambda;\ll)$ with $R(\lambda)$. 

\begin{cor}
    For each 
    $\lambda>\omega$, 
    \begin{equation}\label{E: app_L}
        \lll =\lambda^2\int_0^\infty  e^{-\lambda t}\K_tdt -\lambda\id
    \end{equation}
and $\lll\sra\ll$ on $C^1(\X)$. 
\end{cor}

\subsection{Finite Time-Horizon Approximation}
The current form of the improper integral in \eqref{E: app_L} cannot be directly used for a data-driven approximation. We need to further derive an approximation approach based on finite-horizon observable data. Below, we present a direct truncation modification based on \eqref{E: app_L}.
 
\begin{deff}\label{def: T_t}
    For any $h\in C(\X)$ and $\tau\geq0$, we define $\T_\tau: C(\X)\ra C(\X)$ as
\begin{equation}
    \T_\tau h(x):=\int_0^\tau e^{-\lambda s}\K_{s} h(x) ds. 
\end{equation} 
\end{deff}

It can be shown that for any     $\lambda$, the truncation in Definition \ref{def: T_t} results in an error term $e^{-\lambda \tau} R(\lambda)h(\phi(\tau,x))$  that is uniformly bounded for each observable function $h$, with the dominant term being $e^{-\lambda \tau}$. This indicates that for any arbitrarily large $\lambda$, this truncation does not significantly affect the accuracy compared to the procedure from $\ll_\lambda$ to $\ll$.

For any fixed $\tau> 0$, we  can therefore    use 
\begin{equation}\label{E: approx_t}
  \widetilde{\ll}_\lambda:=  \lambda^2\T_\tau -\lambda\id
\end{equation}
to approximate $\ll$ within a small time-horizon.  We illustrate this approximation with the following simple example.

\begin{eg}
    Consider the simple dynamical system $\dot{\xb}(t)= \xb(t)$ and the observable function $V(x)=x^n$ for any $n\geq 1$. Then, analytically, $\phi(\tau,x)=xe^\tau$ and $\ll V(x) = nx^{n}$. We test the validity of using Eq. \eqref{E: approx_t}. Note that, for sufficiently large $\lambda$, we have 
    \begin{equation*}
        \begin{split}
            &\lambda^2\int_0^\tau  e^{-\lambda s}(\K_sV(x))ds  \\= &\lambda^2 \int_0^\tau e^{-\lambda s}e^{sn}x^n ds =\lambda^2x^n \int_0^\tau e^{-(\lambda-n) s}ds\\
            = & \frac{\lambda^2 x^n}{\lambda -n}(1-e^{-(\lambda -n)\tau})\approx \frac{\lambda^2 x^n}{\lambda-n}
        \end{split}
    \end{equation*}
and 
\begin{equation*}
    \begin{split}
       & \lambda^2\T_\tau V(x)-\lambda V(x) 
       \approx  \frac{\lambda^2 x^n}{\lambda-n}-\lambda x^n 
       \approx nx^n=\ll V(x).
    \end{split}
\end{equation*} 
With high-accuracy evaluation of the integral, we can achieve a reasonably good approximation. \Qed
\end{eg}

\section{Data-Driven Algorithm}\label{sec: alg}
We continue by considering a data-driven algorithm to learn $\widetilde{\ll}_\lambda$, which indirectly approximate $\ll$ as demonstrated in Section \ref{sec: characterize}. 

The idea is to modify the conventional Koopman operator learning scheme  \cite{williams2015data, mauroy2019koopman, meng2023learning}, which involves learning linear operators using a finite-dimensional dictionary of basis functions and representing the image function as a linear combination of the dictionary functions. To be more specific, obtaining a fully discretized version  $\Lb$  of the bounded linear operator $\lambda^2\T_\tau^N-\lambda \id$ based on the training data typically relies on the selection of a discrete dictionary of continuously differentiable observable test functions, denoted by


\begin{equation}\label{E: dict}
    \Zk_N(x):=\left[\zk_0(x),\zk_1(x),\cdots,\zk_{N-1}(x)\right], \;N\in\N.
\end{equation}
Then, the followings should hold: 

1) Let $(\mu_i, \vb_i)_{i=0}^{N-1}$ be the eigenvalues and eigenvectors of $\Lb$. Let $(\rho_i, \varphi_i)_{i=0}^{N-1}$ be the eigenvalues and eigenfunctions of $\ll$. Then, for each $i$, 
\begin{equation}\label{E: approx_eigen}
    \mu_i\approx\rho_i,\;\;\varphi_i(x)\approx \Zk_N(x)\vb_i. 
\end{equation}

2) For any $h\in\operatorname{span}\{\zk_0, \zk_1, \cdots,\zk_{N-1}\}$ such that $h(x)=\Zk_N(x)\mathbf{w}$ for some column vector $\mathbf{w}$, we have that 
\begin{equation}\label{E: approx_Lh}
   \ll h(\cdot)\approx \widetilde{\ll}_\lambda h(\cdot)\approx\lambda^2\T_\tau^Nh(\cdot)-\lambda h(\cdot)\approx \Zk_N(\cdot)(\Lb\mathbf{w}).
\end{equation}
In this section, we modify the existing Koopman learning technique to obtain $\Lb$.

\subsection{Generating Training Data}
For any fixed $\lambda>0$, given a dictionary $\Zk_N$ of the form \eqref{E: dict}, for each $\zk_i\in\Zk_N$ and each $x\in\X$, we consider $\zk_i(x)$ as the features and $\lambda^2\T_\tau\zk_i(x)-\lambda\zk_i(x)=\lambda^2\int_0^\tau e^{-\lambda s} \zk_i(\phi(s,x)) ds-\lambda\zk_i(x)$ as the labels. To compute the integral, we employ numerical quadrature techniques for approximation. This approach inevitably requires discrete-time observations (snapshots), the number of which is denoted by $\Nb$, within the interval $[0,\tau]$ of the flow map $\phi(\tau,x)$. 

To streamline the evaluation process for numerical examples, drawing inspiration from \cite{kang2023data}, for any $\tau$ and $i$, we can assess both the trajectory  and the integral, i.e. the pair $\left(\phi(\tau,x), \lambda^2\int_0^\tau e^{-\lambda s} \zk_i(\phi(s,x)) ds\right)$,  by numerically solving the following  augmented ODE system 
\begin{equation}\label{E: augmented}
    \begin{split}
        &\dot{\xb}(t)=f(\xb(t)), \;\;\xb(0)=x\in\R^n,\\
        &\dot{I}_i(t) = \lambda^2 e^{-\lambda t} \zk_i(\xb),\;\;I_i(0)=0.
    \end{split}
\end{equation}


We summarize the algorithm for generating training data for one time period as in Algorithm \ref{alg: one_time}.

\setcounter{algorithm}{0}
\begin{algorithm}
	\caption{Generating Training Data}\label{alg: one_time} 
	\begin{algorithmic}
		\Require  $f$, $n$, $N$, $\X$, $\Zk_N$, $\tau$, 
  $\lambda$,  and uniformly sampled $\set{x^{(m)}}_{m=0}^{M-1}\subseteq \X$.
  \For{$m$ \textbf{from} 0 \textbf{to} $M-1$}
  \For{$i$ \textbf{from} 0 \textbf{to} $N-1$}
  \State Compute $\zk_i(x^{(m)})$
   \State Compute $\lambda^2\T_\tau\zk_i(x^{(m)})$  using  \eqref{E: augmented} 
 \EndFor
 \State Stack 
 \begin{small}
      \begin{equation*}
      \begin{split}
          \widetilde{\ll}_\lambda\Zk_N(x^{(m)})&=[\lambda^2\T_\tau\zk_0(x^{(m)})-\lambda\zk_0(x^{(m)}),  \cdots,\\
          &\lambda^2\T_\tau\zk_{N-1}(x^{(m)})-\lambda\zk_{N-1}(x^{(m)})]
      \end{split}
 \end{equation*}
 \end{small}
 \EndFor
 \State Stack $\mathfrak{X}, \mathfrak{Y}\in\C^{M\times N }$ such that 
 $\mathfrak{X}=[\Zk_N(x^{(0)}), \Zk_N(x^{(1)}),\cdots, \Zk_N(x^{(M-1)})]^T$
 and
 $\mathfrak{Y}=[\widetilde{\ll}_\lambda\Zk_N(x^{(0)}), \widetilde{\ll}_\lambda\Zk_N(x^{(1)}),\cdots, \widetilde{\ll}_\lambda\Zk_N(x^{(M-1)})]^T$
	\end{algorithmic}
\end{algorithm}

\subsection{Extended Dynamic Mode Decomposition Algorithm}\label{sec: edmd}
After obtaining the training data $(\mathfrak{X}, \mathfrak{Y})$ using Algorithm \ref{alg: one_time}, we can find $\Lb$ by  $\Lb = \operatorname{argmin}_{A\in\C^{N\times N}}\|\mathfrak{Y}-\mathfrak{X}A\|_F$. The $\Lb$ is given in closed-form as $\Lb = \left(\mathfrak{X}^T\mathfrak{X}\right)^\dagger\mathfrak{X}^T\mathfrak{Y}$, 
where $\dagger$ is the pseudo inverse. 
Similar to Extended Dynamic Mode Decomposition (EDMD) \cite{williams2015data} for learning Koopman operators, 
the approximations \eqref{E: approx_eigen} and \eqref{E: approx_Lh} can be guaranteed. In addition, by the universal approximation theorem, all of the function approximations from above should have uniform convergence to the actual quantities. 

It is worth noting that operator learning frameworks, such as autoencoders, which utilize neural networks (NN) as dictionary functions, can reduce human bias in the selection of these functions \cite{deka2022koopman}. Incorporating the  proposed learning scheme with NN falls outside the scope of this paper but will be pursued in future work.

\begin{rem}[The Logarithm Method]\label{rem: log_explain}
    We compare the aforementioned  learning approach of $\ll$  to  those obtained using the benchmark approach as described in \cite{mauroy2019koopman}. Briefly speaking, at a given $s\in[0, \tau]$, that method first obtains a matrix $\mathbf{K}\in\C^{N\times N}$ such that $\mathbf{K} = \left(\mathfrak{X}^T\mathfrak{X}\right)^\dagger\mathfrak{X}^T\widetilde{\mathfrak{Y}}$, where $\widetilde{\mathfrak{Y}}=[\Zk_N(\phi(s,x^{(0)}), \Zk_N(\phi(s,x^{(1)}), \cdots, \Zk_N(\phi(s,x^{(M-1)})]^T.$
    Let $(\mu_i^K, \vb_i^K)_{i=0}^{N-1}$ be the eigenvalues and eigenvectors of $\mathbf{K}$. Let $(\rho_i^K, \varphi_i^K)_{i=0}^{N-1}$ be the eigenvalues and eigenfunctions of $\K_s$.
    Similar to $\eqref{E: approx_eigen}$ and \eqref{E: approx_Lh}, for each $i$, we have $\varphi_i^K(x)\approx \Zk_N(x)\xi_i^K$ and $ \K_s h(\cdot)\approx \Zk_N(\cdot)(\mathbf{K}\mathbf{w})$ for $h(x)=\Zk_N(x)\mathbf{w}$. 

    However, even when $\ll$ can be represented by $(\log{\K_s})/s$, we cannot guarantee that 
    $\frac{\log{\K_s}}{s}h(\cdot)\approx \Zk_N(\cdot)(\frac{\log(\mathbf{K})}{s}\mathbf{w})$, not to mention the case where the above 
    logarithm representation does not hold.  Denoting $\Phi(\cdot)=[\varphi_0^K(\cdot), \varphi_1^K(\cdot),\cdots,\varphi_{N-1}^K(\cdot)]$ and $\Xi=[\xi_0^K, \xi_1^K, \cdots,\xi_{N-1}^K]$, then it is clear that $\Phi(\cdot) = \Zk_N(\cdot)\Xi$. 
    The (possibly complex-valued) rotation matrix $\Xi$ establishes the connection between finite-dimensional eigenfunctions and dictionary functions through data-fitting, ensuring that
    any linear combination within $\Zk_N$ can be equivalently represented using $\Phi$ with a cancellation of the imaginary parts. 
    
This imaginary-part cancellation effect does not generally hold when applying the matrix logarithm. Suppose the imaginary parts account for a significantly large value, the mutual representation of $\Phi$ and $\Zk_N$ does not match in the logarithmic scale. An exception holds unless the chosen dictionary is inherently rotation-free with respect to the true eigenfunctions \cite{zeng2022sampling}, or there is direct access to the data for  $\log(\K_s)$ allowing for direct training of the matrix. However, such conditions contravene our objective of leveraging Koopman data to conversely find the  generator. In comparison, the approach  in this paper presents an elegant method for  approximating $\ll$  regardless of its boundedness. This enables the direct learning of $\Lb$ without  computing the logarithm, thereby avoiding the potential appearance of imaginary parts caused by basis rotation. \Qed
\end{rem}

\section{Case Study}\label{sec: num}
We provide a numerical example to demonstrate the effectiveness of the proposed approach. The research code can be found at \url{https://github.com/Yiming-Meng/Log-Free-Learning-of-Koopman-Generators}. 

Consider the Van der Pol oscillator 
\begin{equation*}
\dot \xb_1(t) = \xb_2(t), \quad 
\dot \xb_2(t) = -\xb_1(t) + (1 - \xb_1^2(t))\xb_2(t), 
\end{equation*} 
with $\xb(0)=x:=[x_1, x_2]$. We assume the system dynamics are unknown to us, and our information is limited to the system dimension, $n=2$, and observations of sampled trajectories. 
To generate training data using Algorithm \ref{alg: one_time}, for simplicity of illustration, we  select $\X=(-1, 1)^2$ and obtain a total of $M=100^2$  uniformly
spaced samples $\set{x^{(m)}}_{m=0}^{M-1}$ in $\X$.  We choose the dictionary as 
\begin{equation}
    \begin{split}
        \Zk_N & = [\zk_{0, 0}, \zk_{1, 0}, \cdots, \zk_{3,0}, \\
        & \zk_{0,1}, \zk_{1,1}, \zk_{2,1}, \cdots, \zk_{i, j}, \cdots,\zk_{3, 2}], \;N=12,
    \end{split}
\end{equation}
 where 
 $\zk_{i, j}(x_1, x_2) = x_1^ix_2^j$ for each $i\in\set{0, 1, 2, 3}$ and $j\in\set{0, 1, 2}$.  We also set  $\tau=1$ and 
$\lambda=10^6$. 
The discrete form $\Lb$ can be obtained according to Section \ref{sec: edmd}.  We apply the learned $\Lb$ to identify the system vector fields, and construct a local Lyapunov function for the unknown system. 

\subsection{System Identification}
The actual vector field is $$f(x) :=[f_1(x), f_2(x)]= [x_2, -x_1 + (1-x_1^2)x_2].$$ As we can analytically establish that $\ll\zk_{1,0}(x) = f_1(x)$ and $\ll\zk_{0,1}(x) = f_2(x)$, we use the  approximation $[\widetilde{\ll}_\lambda \zk_{1,0}, \widetilde{\ll}_\lambda \zk_{0,1}]$ to conversely obtain $f$.

Note that $\zk_{1,0}(x) = \mathcal{Z}_N(x) \mathbf{e}_2$ and $\zk_{0,1}(x) = \mathcal{Z}_N(x) \mathbf{e}_5$, where  each $\mathbf{e}_i$ for $i\in\set{0, 1, \cdots, N}$ is a column unit vector with all components being $0$ except for the $i$-th component, which is $1$.  
To apply $\Lb$, we have that 
$\widetilde{\ll}_\lambda\zk_{1,0}(x) \approx \mathcal{Z}_N(x)(\Lb\mathbf{e}_2)$ and $\widetilde{\ll}_\lambda\zk_{0,1}(x) \approx \mathcal{Z}_N(x)(\Lb\mathbf{e}_5)$. 

In other words, we approximate $\widetilde{\ll}_\lambda \zk_{1,0}$ and $\widetilde{\ll}_\lambda \zk_{0, 1}$ (and hence $f_1$ and $f_2$) using a linear combination of functions within $\Zk_N$. The weights for these approximations are given by $\Lb\mathbf{e}_2$ and $\Lb\mathbf{e}_5$, respectively. We report the corresponding results in Table \ref{tab:f1} and \ref{tab:f2}.

\begin{table}
    \centering
    \begin{tabular}{c|ccc}
        $\zk_{i,j}$ & $j=0$ & $j=1$ & $j=2$ \\\hline
        $i=0$ & $-1.6\times10^{-10}$ & $\mathbf{1+10^{-6}}$  & $-4.8\times10^{-10}$ \\
        $i=1$ & $-1.3\times10^{-5}$  & $-1.9\times10^{-9}$ & $-1.4\times10^{-7}$  \\
        $i=2$ & $-3.6\times10^{-10}$ & $-1.0\times10^{-6}$ & $-6.1\times10^{-10}$  \\
        $i=3$ & $-1.5\times10^{-8}$  & $-2.3\times10^{-9}$ & $-2.9\times10^{-7}$
    \end{tabular}
    \caption{The weights of $\Zk_N$ obtained by $\Lb \mathbf{e}_2$. }
    \label{tab:f1}
\end{table}

\begin{table}
    \centering
    \begin{tabular}{c|ccc}
        $\zk_{i,j}$ & $j=0$ & $j=1$ & $j=2$ \\\hline
        $i=0$ & $-3.8\times10^{-10}$ & $\mathbf{1-10^{-5}}$  & $1.3\times10^{-9}$ \\
        $i=1$ & $\mathbf{-1-10^{-6}}$  & $9.5\times10^{-10}$ & $-2.0\times10^{-6}$  \\
        $i=2$ & $8.1\times10^{-10}$ & $\mathbf{-1-10^{-6}}$ & $-2.0\times10^{-9}$  \\
        $i=3$ & $1\times10^{-8}$  & $-6.9\times10^{-10}$ &$-3.0\times10^{-9}$
    \end{tabular}
    \caption{The weights of $\Zk_N$ obtained by $\Lb \mathbf{e}_5$. }
    \label{tab:f2}
\end{table}

We compare the aforementioned results with those obtained using the benchmark approach as described in \cite{mauroy2019koopman} with the same $M$ and $\Zk_N$. As described in Remark \ref{rem: log_explain}, $\mathbf{K}\in\C^{N\times N}$ can be obtained 
such that $\K_s\zk_{1,0}(x)\approx \Zk_N(x)(\mathbf{K}\mathbf{e}_2)$ and $\K_s\zk_{0,1}(x)\approx \Zk_N(x)(\mathbf{K}\mathbf{e}_5)$. Then, we have the approximation, as claimed by  \cite{mauroy2019koopman}, $f_1(x)\approx \Zk_N(x)(\log(\mathbf{K})\mathbf{e}_2/s)$ and $f_2(x)\approx \Zk_N(x)(\log(\mathbf{K})\mathbf{e}_5/s)$ w.r.t. $|\cdot|_\infty$.   According to \cite[Section VI.A]{mauroy2019koopman}, we set $s=0.5$ to avoid the multi-valued matrix logarithm. We report the weights obtained by taking the real parts of $\log(\mathbf{K})\mathbf{e}_2/s$ and $\log(\mathbf{K})\mathbf{e}_5/s$ in Table \ref{tab:f1_K} and \ref{tab:f2_K}, respectively. Multiple orders of magnitude in accuracy have been established using the proposed method.

It is worth noting that, unlike the experiment in \cite{mauroy2019koopman} where $i,j\in\set{0, 1, 2, 3}$, we deliberately choose different numbers for $i$ and $j$, leading to non-negligible imaginary parts after taking the matrix logarithm of $\mathbf{K}$.  The imaginary parts of the learned weights using the Koopman-logarithm approach are reported in Table \ref{tab:f1_K_im} and \ref{tab:f2_K_im}.  The presence of non-negligible imaginary parts indicates that the underlying system is sensitive to the selection of dictionary functions when employing the Koopman-logarithm approach,  an effect that is significant and cannot be overlooked when aiming to minimize human intervention in identifying unknown systems in practice.
Furthermore, the original experiment in \cite{mauroy2019koopman} used $M=300^2$ samples for a data fitting, and when $M$ is reduced to $100^2$, the quality of the operator learning deteriorates.

\begin{table}
    \centering
    \begin{tabular}{c|ccc}
        $\zk_{i,j}$ & $j=0$ & $j=1$ & $j=2$ \\\hline
        $i=0$ & $-9.5\times10^{-16}$ & $\mathbf{1-3.0\times 10^{-4}}$  & $-4.8\times10^{-10}$ \\
        $i=1$ & $\mathbf{-4.1\times10^{-3}}$  & $-3\times10^{-15}$ & $\mathbf{-1.8\times10^{-2}}$   \\
        $i=2$ & $-1.7\times10^{-16}$ & $\mathbf{-6.4\times10^{-3}}$ &  $-7.8\times10^{-17}$\\
        $i=3$ & $\mathbf{-6.0\times10^{-3}}$  & $-6.7\times10^{-15}$ & $\mathbf{-3.7\times10^{-2}}$ 
    \end{tabular}
    \caption{The weights of $\Zk_N$ obtained by $\log(\mathbf{K})\mathbf{e}_2/s$. }
    \label{tab:f1_K}
\end{table}

\begin{table}
    \centering
    \begin{tabular}{c|ccc}
        $\zk_{i,j}$ & $j=0$ & $j=1$ & $j=2$ \\\hline
        $i=0$ & $-8.6\times10^{-16}$ & $\mathbf{1+6.5\times 10^{-3}}$  & $4.1\times10^{-16}$ \\
        $i=1$ & $-\mathbf{1+3.0\times10^{-2}}$  & $-1.6\times10^{-15}$ & $\mathbf{-1.5\times10^{-1}}$   \\
        $i=2$ & $1.0\times10^{-15}$ & $-\mathbf{1-3.8\times 10^{-2}}$ &  $-2.0\times10^{-15}$\\
        $i=3$ & $\mathbf{-5.8\times10^{-2}}$  & $1.8\times10^{-15}$ & $\mathbf{3.3\times10^{-1}}$ 
    \end{tabular}
    \caption{The weights of $\Zk_N$ obtained by $\log(\mathbf{K})\mathbf{e}_5/s$ }
    \label{tab:f2_K}
\end{table}

\begin{table} 
    \centering
    \begin{tabular}{c|ccc}
        $\zk_{i,j}$ & $j=0$ & $j=1$ & $j=2$ \\\hline
        $i=0$ & $-1.8\times10^{-17}$ & $\mathbf{2.0\times 10^{-3}}$  & $-2.2\times10^{-18}$ \\
        $i=1$ & $\mathbf{5.3\times10^{-3}}$  & $6.9\times10^{-17}$ & $\mathbf{-2.3\times10^{-2}}$   \\
        $i=2$ & $1.6\times10^{-17}$ & $\mathbf{-8.4\times10^{-3}}$ &  $-3.0\times10^{-17}$\\
        $i=3$ & $\mathbf{-8.4\times10^{-3}}$  & $-1.1\times10^{-16}$ & $\mathbf{5.0\times10^{-2}}$ 
    \end{tabular}
    \caption{The imaginary parts obtained by $\log(\mathbf{K})\mathbf{e}_2/s$. }
    \label{tab:f1_K_im}
\end{table}

\begin{table} 
    \centering
    \begin{tabular}{c|ccc}
        $\zk_{i,j}$ & $j=0$ & $j=1$ & $j=2$ \\\hline
        $i=0$ & $-8.0\times10^{-17}$ & $\mathbf{8.8\times 10^{-3}}$  & $-9.5\times10^{-18}$ \\
        $i=1$ & $\mathbf{2.3\times10^{-2}}$  & $3.0\times10^{-16}$ & $\mathbf{-1.0\times10^{-1}}$   \\
        $i=2$ & $6.9\times10^{-17}$ & $-\mathbf{-3.7\times 10^{-2}}$ &  $-1.3\times10^{-16}$\\
        $i=3$ & $\mathbf{3.7\times10^{-2}}$  & $-4.9\times10^{-16}$ & $\mathbf{2.2\times10^{-1}}$ 
    \end{tabular}
    \caption{The imaginary parts obtained by obtained by $\log(\mathbf{K})\mathbf{e}_5/s$.}
    \label{tab:f2_K_im}
\end{table}

\subsection{Stability Prediction for the Reversed Dynamics}
\begin{figure}[!t]
\centerline{\includegraphics[scale = 0.47
]{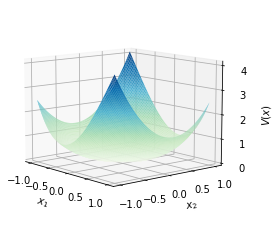} 
\includegraphics[scale = 0.47
]{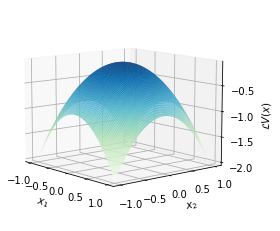}}
\caption{Data-driven Lyapunov function $V$ and its Lie derivative $\ll V$ w.r.t. the predicted system dynamics. }
\label{fig: V}
\end{figure}

Observing the identified system dynamics, we anticipate that the time-reversed system is (locally) asymptotically stable w.r.t. the origin. We use the learned $\Lb$ to construct polynomial Lyapunov functions based on the Lyapunov equation $\ll V(x)=|x|^2$, where the sign on the r.h.s. is reversed due to the dynamics being reversed. To use $\Lb$ and the library functions from $\Zk_N$, we define the weights as $\theta:=[\theta_0,\theta_1, \cdots,\theta_{N-1}]^T$ and seek a $\theta$ such that  $V(x;\theta):=\Zk_N(x)\theta$ , with the objective of minimizing $|\Zk_N(x)\Lb \theta-x_1^2-x_2^2|$. Ignoring the small terms of the magnitude $10^{-9}$, the constructed Lyapunov function is given by $$V(x;\theta) = 1.39x_1^2 - 1.56 x_1x_2 + 1.16x_2^2 + 0.74x_1^2x_2^2. $$ Additionally, $\ll V(x;\theta)$ is approximated by $\Zk_N(x)\Lb\theta$, with the maximal value observed at  $\Zk_N(0)\Lb\theta = -0.068$, which indicates that stability predition using the data-driven Lyapunov function is verified to be valid. This indicates that the stability prediction, as inferred using the data-driven Lyapunov function, is verified to be valid. The visualization can be found in Fig \ref{fig: V}. 

\begin{rem}
    The idea illustrated above can be expanded to cover a larger region of interest. Specifically, a Zubov equation, instead of a Lyapunov equation,  can be solved within the Koopman learning framework using the same dataset \cite{meng2023learning} as proposed in Algorithm \ref{alg: one_time}. The solution obtained can potentially serve as a Lyapunov function. Since it cannot guarantee the properties of the learned  function's derivatives, to confirm it as a true Lyapunov function, the data can be reused as in Algorithm \ref{alg: one_time} to verify its Lie derivative. \Qed
\end{rem}

\section{Conclusion}\label{sec: conclusion}
In this paper, we propose a logarithm-free Koopman operator-based learning framework for the infinitesimal generator, demonstrating both theoretical and numerical improvements over the method proposed in \cite{mauroy2019koopman}. In particular, for more general cases where the generator is unbounded and, consequently, the logarithm of the Koopman operator cannot be used for representation, we draw upon the rich literature to propose an approximation (Eq. \eqref{E: approx_t}) based on Yosida's approximation. A convergence result, along with the convergence rate, is proved in this paper to guide users in tuning the parameters. A numerical example with application in system identification is provided in comparison with the experiment in \cite{mauroy2019koopman}  demonstrating 
the learning accuracy. Unlike the experiment in \cite{mauroy2019koopman}, where learning accuracy is sensitive to the choice of dictionary functions, the method presented in this paper shows significant improvement in this regard. In applications where automatic computational approaches surpass human computability, such as in constructing Lyapunov-like functions using the Lie derivative, the proposed logarithm-free method holds more promise. We will pursue future efforts to provide more analysis on the sampling rate, numerical simulations, and real-world applications using real data.

\bibliographystyle{ieeetr}        
\bibliography{TAC}

\end{document}